# A FINITE BASIS THEOREM FOR RESIDUALLY FINITE, CONGRUENCE MEET-SEMIDISTRIBUTIVE VARIETIES

ROSS WILLARD

ABSTRACT. We derive a Mal'cev condition for congruence meet-semidistributivity and then use it to prove two theorems. **Theorem A**: if a variety in a finite language is congruence meet-semidistributive and residually less than some finite cardinal, then it is finitely based. **Theorem B**: there is an algorithm which, given $m < \omega$ and a finite algebra in a finite language, determines whether the variety generated by the algebra is congruence meet-semidistributive and residually less than $m$.

## 1. INTRODUCTION

We consider finite algebras and the varieties they generate. A problem which is currently of interest to general algebraists is the following speculation, dating from the mid 1970s (see Conjecture 1 in R. Park's Ph.D. thesis [21]) and attributed, perhaps erroneously, to B. Jónsson in [13]:

**Problem 1.** Suppose $\mathbf{A}$ is a finite algebra in a finite language, and suppose there exists $m < \omega$ such that every subdirectly irreducible member of $\mathbf{HSP}(\mathbf{A})$ has cardinality less than $m$. Must the equational theory of $\mathbf{A}$ be finitely based?

A variety for which every subdirectly irreducible member has cardinality less than $m$ for some $m < \omega$ is said to *have a finite residual bound*. For example, every congruence distributive variety generated by a finite algebra has a finite residual bound, by an old result of Jónsson [6, Corollary 2.5] (see also [4, Theorem 2.5]). This fact played a key role in the proof of K. Baker's celebrated finite basis theorem [1] of 1972, which gives a positive answer to Problem 1 in case $\mathbf{HSP}(\mathbf{A})$ is congruence distributive and was surely part of the motivation for the problem. The problem received partial confirmation in 1985 through R. McKenzie's extension of Baker's theorem to congruence modular varieties [14]. In this paper we give further confirmation by extending Baker's theorem in another direction.

A second problem currently of interest is the so-called RS problem. This asks for a characterization of those finite algebras $\mathbf{A}$ for which $\mathbf{HSP}(\mathbf{A})$ is residually

1991 *Mathematics Subject Classification.* Primary 08B05; Secondary 03C05.
*Key words and phrases.* equational theory, finitely based, residually finite, variety, congruence meet-semidistributive.
The financial support of the NSERC of Canada is gratefully acknowledged.

small. The characterization cannot be decidable: McKenzie has proved [18] that the properties "**HSP**(**A**) is residually large" and "**HSP**(**A**) has a finite residual bound" are recursively inseparable properties of **A**. However, it is expected that a full solution to the RS problem will yield as a byproduct a proof that the former property is recursively enumerable. In fact, we conjecture that both properties are r.e., and in the present paper confirm that the latter property is r.e. when restricted to the class of algebras covered by our finite basis theorem.

A variety is *congruence meet-semidistributive* if the congruence lattice of each member satisfies the meet-semidistributive law $x \wedge y = x \wedge z \to x \wedge y = x \wedge (y \vee z)$. Congruence meet-semidistributive varieties include congruence distributive varieties as well as varieties of semilattice-based algebras. G. Czédli proved in 1982 that congruence meet-semidistributivity is a weak Mal'cev property of varieties. Recently K. Kearnes and Á. Szendrei [10] and P. Lipparini [12] proved that it is in fact a Mal'cev property. In the present paper we give a direct proof of this fact, and then show that the corresponding Mal'cev condition provides a natural setting for the combinatorial arguments in Baker's original proof of his theorem. Thus we are able to prove a finite basis theorem for finite algebras **A** in a finite language for which **HSP**(**A**) is congruence meet-semidistributive and has a finite residual bound. We then use these combinatorial arguments to prove that the property "**HSP**(**A**) has a finite residual bound" is an r.e. property of **A** when restricted to algebras **A** for which **HSP**(**A**) is congruence meet-semidistributive. This complements McKenzie's proof [19] that the property "**HSP**(**A**) is residually large" is r.e. when restricted to these same algebras.

Finally, let us mention the manuscript [11], in which similar arguments are used to prove that every residually finite, congruence meet-semidistributive variety in a finite language has a finite residual bound. This explains the title of our paper.

This paper grew out of our efforts [23], later with G. McNulty [20], to prove finite basis theorems for certain finite semilattice-based algebras related to McKenzie's refutation of the RS conjecture. In particular, the finite basis theorem in this paper makes our paper [23] obsolete and thus gives (with [16]) yet another route to McKenzie's negative solution to Tarski's finite basis problem [17].

I wish to thank the Fields Institute for Research in Mathematical Sciences in Toronto, Canada, which provided the ideal setting for the discovery of the results in this paper. I also thank Kirby Baker, Joel Berman, Gábor Czédli, Keith Kearnes, Paolo Lipparini, and especially George McNulty for their helpful discussions and comments.

2. Characterizing congruence meet-semidistributivity

Given an algebra $\mathbf{A}$ and congruences $\alpha, \beta, \gamma$ of $\mathbf{A}$, define

$$\begin{aligned}
\beta_0 &= \beta, & \gamma_0 &= \gamma, \\
\beta_{n+1} &= \beta \vee (\alpha \cap \gamma_n), & \gamma_{n+1} &= \gamma \vee (\alpha \cap \beta_n), \\
\beta_\omega &= \bigcup_{n<\omega} \beta_n, & \gamma_\omega &= \bigcup_{n<\omega} \gamma_n.
\end{aligned}$$

$\beta_\omega$ and $\gamma_\omega$ are the least congruences of $\mathbf{A}$ extending $\beta$ and $\gamma$ respectively and satisfying $\alpha \cap \beta_\omega = \alpha \cap \gamma_\omega$. Thus if $\mathbf{A}$ is congruence meet-semidistributive then $\alpha \cap (\beta \vee \gamma) \leq \beta_\omega$. Czédli [3, Claim 3] showed that congruence meet-semidistributivity in an algebra is equivalent to the universally quantified condition $\alpha \cap (\beta \vee \gamma) \leq \beta_\omega$, which in turn is equivalent to the infinite conjunction of the conditions

$$\alpha \cap (\underbrace{\beta \circ \gamma \circ \beta \circ \cdots}_{k-1 \text{ o's}}) \subseteq \beta_\omega, \qquad k \geq 2.$$

For fixed $k \geq 2$, let $\mathcal{C}_k$ be the displayed inclusion. Each $\mathcal{C}_k$ can be characterized at the level of varieties by a Mal'cev condition in the usual way. We shall do this for $\mathcal{C}_2$.

Our characterization is developed using the notion of coloured ordered (finite, rooted) trees. (A tree is *ordered* if for each node $p$ there is an assigned linear ordering of the children of $p$. A tree is *coloured* if it is accompanied by a function $\chi$ whose domain is the set of nodes.) In the trees that we shall consider, the colour-set (codomain of $\chi$) shall be $\{b, g\}$ where informally $b$ and $g$ represent the colours blue and green, or alternatively, the congruences $\beta$ and $\gamma$ respectively in the construction of $\beta_\omega$. Furthermore, we assume that the colour of a child is always opposite to the colour of its parent.

Let $T$ be such a tree. A variety $\mathcal{V}$ satisfies the condition $\mathcal{C}_2(T)$ if there exist two indexed families $\{s_p(x,y,z) : p \in T\}$ and $\{t_p(x,y,z) : p \in T\}$ of 3-ary terms in the language of $\mathcal{V}$ such that the following equations are identically true in $\mathcal{V}$:

1. $s_0(x,y,z) \approx x$ and $t_0(x,y,z) \approx z$ (where 0 is the root of $T$).
2. $s_p(x,y,x) \approx t_p(x,y,x)$ (for all $p \in T$).
3. (If $\chi(p) = b$ and $r$ is the first child of $p$): $s_p(x,x,y) \approx s_r(x,x,y)$.
4. (If $\chi(p) = g$ and $r$ is the first child of $p$): $s_p(x,y,y) \approx s_r(x,y,y)$.
5. (If $\chi(p) = b$ and $r$ is the last child of $p$): $t_p(x,x,y) \approx t_r(x,x,y)$.
6. (If $\chi(p) = g$ and $r$ is the last child of $p$): $t_p(x,y,y) \approx t_r(x,y,y)$.
7. (If $p, q$ are consecutive children, in that order, of some node labeled by $b$): $t_p(x,x,y) \approx s_q(x,x,y)$.
8. (If $p, q$ are consecutive children, in that order, of some node labeled by $g$): $t_p(x,y,y) \approx s_q(x,y,y)$.
9. (If $p$ is a leaf labeled by $b$): $s_p(x,x,y) \approx t_p(x,x,y)$.
10. (If $p$ is a leaf labeled by $g$): $s_p(x,y,y) \approx t_p(x,y,y)$.

For example, let $T$ be the tree pictured below (0 is the root and is coloured by $\mathsf{b}$; the children of 0 are ordered left to right; the colours of the vertices other than 0 alternate along each branch.) The identities of $\mathcal{C}_2(T)$ (suppressing $s_0$ and $t_0$) are

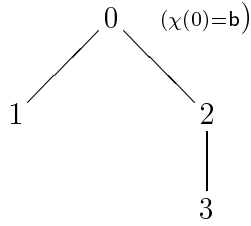

$$\begin{aligned}
s_p(x,y,x) &\approx t_p(x,y,x), \quad p=1,2,3, \\
s_1(x,x,y) &\approx x, \\
s_2(x,y,y) &\approx s_3(x,y,y), \\
t_2(x,x,y) &\approx y, \\
t_2(x,y,y) &\approx t_3(x,y,y), \\
t_1(x,x,y) &\approx s_2(x,x,y), \\
s_3(x,x,y) &\approx t_3(x,x,y), \\
s_1(x,y,y) &\approx t_1(x,y,y).
\end{aligned}$$

The variety of semilattices satisfies $\mathcal{C}_2(T)$ as witnessed by the terms

$$\begin{aligned}
s_1 &= xy, & t_1 &= xyz, \\
s_2 &= xz, & t_2 &= z, \\
s_3 &= xyz, & t_3 &= yz.
\end{aligned}$$

**Theorem 2.1.** *For a variety $\mathcal{V}$, the following are equivalent:*

1. *$\mathcal{V}$ is congruence meet-semidistributive.*
2. *$\mathcal{V}$ satisfies $\mathcal{C}_2$.*
3. *The 3-generated free algebra in $\mathcal{V}$ is congruence meet-semidistributive.*
4. *The 3-generated free algebra in $\mathcal{V}$ satisfies $\mathcal{C}_2$.*
5. *$\mathcal{V}$ satisfies $\mathcal{C}_2(T)$ for some (finite rooted) coloured ordered tree $T$ as above.*
6. *There exists a finite family $\{\langle s_p(x,y,z), t_p(x,y,z)\rangle : p \in T\}$ of pairs of ternary terms such that*

$$\mathcal{V} \models s_p(x,y,x) \approx t_p(x,y,x) \quad (p \in T),$$

$$\mathcal{V} \models \forall xy \left( x = y \leftrightarrow \bigwedge_{p \in T} [s_p(x,x,y) = t_p(x,x,y) \leftrightarrow s_p(x,y,y) = t_p(x,y,y)] \right).$$

7. *For all $\mathbf{A} \in \mathcal{V}$ and $a_0, a_1, \ldots, a_n \in A$, if $a_0 \neq a_n$ then there exists $i < n$ such that $\mathrm{Cg}(a_0, a_n) \cap \mathrm{Cg}(a_i, a_{i+1}) \neq 0_A$.*

**Remark.** The equivalence of items 1 through 4 was first proved by Kearnes and Szendrei [10] and by Lipparini [12]. Earlier, D. Hobby and McKenzie had proved it for locally finite varieties [5].

*Proof of Theorem 2.1.* $(1 \Rightarrow 3)$ and $(2 \Rightarrow 4)$ are trivial, while $(1 \Rightarrow 2)$ and $(3 \Rightarrow 4)$ follow from Czédli's argument.

($4 \Rightarrow 5$). Let $\mathbb{F}_\mathcal{V}(x,y,z)$ be the $\mathcal{V}$-free algebra on the generators $x, y, z$ and assume $\mathbb{F}_\mathcal{V}(x,y,z)$ satisfies $\mathcal{C}_2$. Define

$$\begin{aligned} \alpha &= \mathrm{Cg}(x,z), \\ \beta &= \mathrm{Cg}(x,y), \\ \gamma &= \mathrm{Cg}(y,z). \end{aligned}$$

Thus $(x,z) \in \alpha \cap (\beta \circ \gamma)$. As $\mathbb{F}_\mathcal{V}(x,y,z)$ satisfies $\mathcal{C}_2$ we have $(x,z) \in \beta_\omega$ and thus $(x,z) \in \alpha \cap \beta_n$ for some $n < \omega$. We adopt the usual convention that ternary terms $s(x,y,z)$ name elements $s$ of $\mathbb{F}_\mathcal{V}(x,y,z)$, so that $s = t$ in $\mathbb{F}_\mathcal{V}(x,y,z)$ if and only if $\mathcal{V} \models s(x,y,z) \approx t(x,y,z)$. The idea is to now show by induction on $k$ that

1. If $(s,t) \in \alpha \cap \beta_k$ then there exist a coloured ordered tree $T$ of height $k$ and with root 0, and two indexed families $\{s_p(x,y,z) : p \in T\}$ and $\{t_p(x,y,z) : p \in T\}$ of 3-ary terms in the language of $\mathcal{V}$, such that
    (a) $\chi(0) = \mathsf{b}$.
    (b) The equations in the definition of $\mathcal{C}_2(T)$, with the exception of $s_0(x,y,z) \approx x$ and $t_0(x,y,z) \approx z$, are identically true in $\mathcal{V}$.
    (c) $s_0 = s$ and $t_0 = t$.
2. Similarly, if $(s,t) \in \alpha \cap \gamma_k$ then there exist $T$ and the two families of terms satisfying the same requirements as in the first item except that $\chi(0) = \mathsf{g}$.

When $k = 0$ we can choose $T$ to be the 1-element tree. If $(s,t) \in \alpha \cap \beta_{k+1}$ then since $(s,t) \in \beta \vee (\alpha \cap \gamma_k)$ there exist $s_1, t_1, \ldots, s_m, t_m \in \mathbb{F}_\mathcal{V}(x,y,z)$ so that $s \overset{\beta}{\equiv} s_1 \overset{\alpha \cap \gamma_k}{\equiv} t_1 \overset{\beta}{\equiv} s_2 \overset{\alpha \cap \gamma_k}{\equiv} \cdots \overset{\beta}{\equiv} s_m \overset{\alpha \cap \gamma_k}{\equiv} t_m \overset{\beta}{\equiv} t$. For $i = 1, \ldots, m$ let $T_i$ be a coloured ordered tree witnessing $(s_i, t_i) \in \alpha \cap \gamma_k$ as in item 2; rename the root of $T_i$ by $i$ and assume that $T_1, \ldots, T_m$ and $\{0\}$ are pairwise disjoint. Define $T = \{0\} \cup T_1 \cup \cdots \cup T_m$ with root 0 and $\chi(0) = \mathsf{b}$ and so that the children of 0 are $1, \ldots, m$ in that order. This works; a similar argument works if $(s,t) \in \alpha \cap \gamma_{k+1}$.

($5 \Rightarrow 6$). Let $T$ be a coloured ordered tree such that $\mathcal{V}$ satisfies $\mathcal{C}_2(T)$, witnessed by $\{s_p(x,y,z) : p \in T\}$ and $\{t_p(x,y,z) : p \in T\}$. We shall show that the family $\{\langle s_p(x,y,z), t_p(x,y,z)\rangle : p \in T\}$ satisfies item 6.

Let $\mathbf{A} \in \mathcal{V}$ and $a, b \in A$ be given. It remains to show that $a = b$ under the assumption that $s_p(a,a,b) = t_p(a,a,b) \leftrightarrow s_p(a,b,b) = t_p(a,b,b)$ for all $p \in T$. We argue by induction on $p$ (starting with the leaves) that $s_p(a,a,b) = t_p(a,a,b)$ and $s_p(a,b,b) = t_p(a,b,b)$. If $p$ is a leaf, then one of these equations is guaranteed by $\mathcal{C}_2(T)$, so the other follows by the assumption. Next, assume that $p$ is not a leaf and the claim is true of all the children of $p$. Suppose for the sake of argument that $\chi(p) = \mathsf{b}$. Let $r_1, \ldots, r_\ell$ be the children of $p$ in increasing order. The identities of $\mathcal{C}_2(T)$ and the inductive assumption imply

$$s_p(a,a,b) = s_{r_1}(a,a,b) = t_{r_1}(a,a,b) = s_{r_2}(a,a,b) = \cdots = t_{r_\ell}(a,a,b) = t_p(a,a,b),$$

establishing one of the required equations. The other equation then follows by the assumption. A similar argument works if $\chi(p) = \mathsf{g}$. Thus the equations hold for all $p \in T$. When applied to the root, the first two identities of $\mathcal{C}_2(T)$ yield $a = b$.

($6 \Rightarrow 7$). Given a family $\{\langle s_p(x,y,z), t_p(x,y,z)\rangle : p \in T\}$ of pairs of ternary terms satisfying item 6, let $\mathbf{A} \in \mathcal{V}$ and $a_0, a_1, \ldots, a_n \in A$ be given with $a_0 \neq a_n$ and put $a = a_0$ and $b = a_n$. By the second condition in item 4, there exists $p \in T$ such that $\neg[s_p(a,a,b) = t_p(a,a,b) \leftrightarrow s_p(a,b,b) = t_p(a,b,b)]$. Suppose for concreteness that $s_p(a,a,b) = t_p(a,a,b)$ while $s_p(a,b,b) \neq t_p(a,b,b)$. Choose $i < n$ such that $s_p(a,a_i,b) = t_p(a,a_i,b)$ while $s_p(a,a_{i+1},b) \neq t_p(a,a_{i+1},b)$. Let $c = s_p(a,a_{i+1},b)$, $d = t_p(a,a_{i+1},b)$, $u = s_p(a,a_i,b)$, and $v = s_p(a,a_{i+1},a)$; note that $c \neq d$. Define $f_1(x) = s_p(a,x,b)$, $f_2(x) = t_p(a,x,b)$, $g_1(x) = s_p(a,a_{i+1},x)$, and $g_2(x) = t_p(a,a_{i+1},x)$. Then $f_1, f_2, g_1, g_2 \in \operatorname{Pol}_1 \mathbf{A}$ and

$$\begin{aligned}
\{f_1(a_i), f_1(a_{i+1})\} &= \{c, u\}, \\
\{f_2(a_i), f_2(a_{i+1})\} &= \{u, d\} \quad \text{by the choice of } i, \\
\{g_1(a), g_1(b)\} &= \{c, v\}, \\
\{g_2(a), g_2(b)\} &= \{v, d\} \quad \text{using the identity } s_p(x,y,x) \approx t_p(x,y,x).
\end{aligned}$$

Thus $(c,d) \in \operatorname{Cg}(a_i, a_{i+1}) \cap \operatorname{Cg}(a,b)$.

($7 \Rightarrow 1$). Assume that $\mathcal{V}$ satisfies the condition in item 7. To prove that $\mathcal{V}$ is congruence meet-semidistributive, it suffices to show that if $\mathbf{A} \in \mathcal{V}$ and $\alpha, \beta, \gamma \in \operatorname{Con} \mathbf{A}$ are such that $\alpha \cap \beta = \alpha \cap \gamma = 0_A$, then $\alpha \cap (\beta \vee \gamma) = 0_A$. Suppose instead that $\alpha \cap (\beta \vee \gamma) \neq 0_A$. Thus we may choose $a = a_0, a_1, \ldots, a_n = b$ in $A$ with $a \neq b$ and $(a,b) \in \alpha$ and $(a_i, a_{i+1}) \in \beta \cup \gamma$ for each $i < n$. By item 7, there must exist $i < n$ such that $\operatorname{Cg}(a,b) \cap \operatorname{Cg}(a_i, a_{i+1}) \neq 0_A$. But then $\alpha \cap \beta \neq 0_A$ or $\alpha \cap \gamma \neq 0_A$. □

## 3. The Baker-style argument

Suppose $\mathbf{A}$ is an algebra. Following [1], by a *basic translation* of $\mathbf{A}$ we mean any unary polynomial of the form $F(a_1, \ldots, a_{i-1}, x, a_{i+1}, \ldots, a_n)$ where $F$ is an $n$-ary fundamental operation of $\mathbf{A}$, $1 \leq i \leq n$, and the $a_j$'s are any elements of $A$. A *$k$-translation* of $\mathbf{A}$ is a unary polynomial of $\mathbf{A}$ which can be expressed as the composition of $k$ or fewer basic translations. In particular, the identity map $\operatorname{id}_A$ is the unique 0-translation of $\mathbf{A}$.

$A^{(2)}$ denotes the set of all 2-element subsets of $A$. If $\{a,b\}, \{c,d\} \in A^{(2)}$ and $k < \omega$, then we write $\{a,b\} \to_k \{c,d\}$ to mean that there exists a $k$-translation $f$ such that $\{f(a), f(b)\} = \{c,d\}$. Similarly, if $k, n < \omega$ then we define $\{a,b\} \Rightarrow_{k,n} \{c,d\}$ to mean that there exists a sequence $c = c_0, c_1, \ldots, c_n = d$ such that for each $i < n$, either $c_i = c_{i+1}$ or $\{a,b\} \to_k \{c_i, c_{i+1}\}$. The notation $\{a,b\} \Rightarrow_k \{c,d\}$ means $\{a,b\} \Rightarrow_{k,n} \{c,d\}$ for some $n < \omega$. Note that by Mal'cev's description of principal congruences, if $\{a,b\}, \{c,d\} \in A^{(2)}$ then $(c,d) \in \operatorname{Cg}^{\mathbf{A}}(a,b)$ if and only if $\{a,b\} \Rightarrow_k \{c,d\}$ for some $k < \omega$. Moreover, if the language of **A** is finite then for all $k, n < \omega$ there is a first-order formula $\pi_{k,n}(x, y, z, w)$ (a principal congruence formula in the sense of [22]) which defines the relation $\{x, y\} \Rightarrow_{k,n} \{z, w\}$ in all algebras of the same type as **A**.

The relations $\to_k$ and $\Rightarrow_{k,n}$ have the following properties.

1. $\to_k$ and $\Rightarrow_{k,1}$ mean the same thing.
2. If $\{a, b\} \Rightarrow_{k,m} \{c, d\} \Rightarrow_{\ell,n} \{r, s\}$, then $\{a, b\} \Rightarrow_{k+\ell,mn} \{r, s\}$. In other words, compositions of $\Rightarrow_{x,y}$ are additive in $x$ and multiplicative in $y$.
3. If $\{a, b\} \to_{k+\ell} \{c, d\}$, then there exist $u, v$ such that $\{a, b\} \to_k \{u, v\} \to_\ell \{c, d\}$.

The next two lemmas are patterned after [1, Lemmas 5.5 and 5.3].

**Lemma 3.1** (Single-sequence lemma). *Suppose $\mathcal{V}$ is a class of algebras for which there exist terms $s_p(x, y, z), t_p(x, y, z)$ witnessing Theorem 2.1(6). Suppose moreover that each $s_p, t_p$ is a fundamental operation symbol in the language of $\mathcal{V}$. Then the following is true: if $\mathbf{A} \in \mathcal{V}$ and $a = a_0, a_1, \ldots, a_n = b$ is a sequence in $A$ with $a \neq b$, then there exist $\{c, d\} \in A^{(2)}$ and $i < n$ such that $\{a_i, a_{i+1}\} \Rightarrow_{1,2} \{c, d\}$ and $\{a, b\} \Rightarrow_{1,2} \{c, d\}$.*

*Proof.* The proof of Theorem 2.1(6 $\Rightarrow$ 7) yielded $\{c, d\} \in A^{(2)}$ and basic translations $f_1, f_2, g_1, g_2$ of **A** which witness $\{a_i, a_{i+1}\} \Rightarrow_{1,2} \{c, d\}$ and $\{a, b\} \Rightarrow_{1,2} \{c, d\}$. □

**Lemma 3.2** (Multi-sequence lemma). *With the same assumptions as before, the following is true: if $\mathbf{A} \in \mathcal{V}$, $\{a, b\} \in A^{(2)}$, and $S_1, \ldots, S_N$ are sequences from $a$ to $b$, where $S_i = (a_0^i, a_1^i, \ldots, a_{\lambda(i)}^i)$ with $a_0^i = a$ and $a_{\lambda(i)}^i = b$ for each $i = 1, \ldots, N$, then there exist $\{u, v\} \in A^{(2)}$ and, for each $i$, a 'key link' $\{a_{\sigma(i)}^i, a_{\sigma(i)+1}^i\}$ of distinct adjacent elements of $S_i$, which we shall rename $\{a^i, b^i\}$, such that $\{a^i, b^i\} \Rightarrow_{N,2^N} \{u, v\}$ for each $i$ and $\{a, b\} \Rightarrow_{N,2^N} \{u, v\}$.*

*Proof.* The proof is virtually the same as Baker's proof of his original multisequence lemma [1, Lemma 5.3]. Argue by induction on $N$. If $N = 1$ then the claim is Lemma 3.1. If $N > 1$, apply the claim to the sequences $S_1, \ldots, S_{N-1}$ to get $\{c, d\} \in A^{(2)}$ and key links $\{a^1, b^1\}, \ldots, \{a^{N-1}, b^{N-1}\}$ such that

$$\{a^i, b^i\} \Rightarrow_{N-1, 2^{N-1}} \{c, d\}, \quad i = 1, \ldots, N-1,$$
$$\{a, b\} \Rightarrow_{N-1, 2^{N-1}} \{c, d\}.$$

Choose distinct $c = c_0, c_1, \ldots, c_m = d$ so that $\{a, b\} \to_{N-1} \{c_j, c_{j+1}\}$ for all $j < m$. Choose $(N-1)$-translations $f_0, \ldots, f_{m-1}$ so that $\{c_j, c_{j+1}\} = \{f_j(a), f_j(b)\}$ for all $j < m$. For each $j < m$ define a sequence $T_j$ from $c_j$ to $c_{j+1}$ by applying $f_j$ to $S_N$ or its reverse; that is,

$$T_j = \begin{cases} (f_j(a_0^N), f_j(a_1^N), \ldots, f_j(a_{\lambda(N)}^N)) & \text{if } (f_j(a), f_j(b)) = (c_j, c_{j+1}), \\ (f_j(a_{\lambda(N)}^N), \ldots, f_j(a_1^N), f_j(a_0^N)) & \text{if } (f_j(a), f_j(b)) = (c_{j+1}, c_j). \end{cases}$$

Let $T$ be the sequence from $c$ to $d$ formed by concatenating $T_0, \ldots, T_{m-1}$. By Lemma 3.1 there must exist $\{u, v\} \in A^{(2)}$ and a link $\{r, s\}$ in $T$ such that $r \neq s$, $\{r, s\} \Rightarrow_{1,2} \{u, v\}$ and $\{c, d\} \Rightarrow_{1,2} \{u, v\}$. By construction, the link $\{r, s\}$ must be of the form $\{f_j(a_k^N), f_j(a_{k+1}^N)\}$ for some $j < m$ and $k < \lambda(N)$. Rename $\{a_k^N, a_{k+1}^N\}$ as $\{a^N, b^N\}$, our chosen key link of $S_N$. Then

$$\begin{aligned}
\{a^i, b^i\} &\Rightarrow_{N-1, 2^{N-1}} \{c, d\} \Rightarrow_{1,2} \{u, v\}, \quad i = 1, \ldots, N-1, \\
\{a, b\} &\Rightarrow_{N-1, 2^{N-1}} \{c, d\} \Rightarrow_{1,2} \{u, v\}, \\
\{a^N, b^N\} &\to_{N-1} \{r, s\} \Rightarrow_{1,2} \{u, v\}.
\end{aligned}$$

Thus $\{a^i, b^i\}$ $(1 \leq i \leq N)$ and $\{u, v\}$ have the desired properties. □

**Corollary 3.3.** *With the same hypotheses as in Lemma 3.1, suppose $\mathbf{A} \in \mathcal{V}$, $\{a_1, b_1\}, \ldots, \{a_N, b_N\}, \{u, v\} \in A^{(2)}$ and $n > 0$ are such that $\{a_i, b_i\} \Rightarrow_n \{u, v\}$ for all $i = 1, \ldots, N$. Then there exist $\{r_i, s_i\} \in A^{(2)}$ $(1 \leq i \leq N)$ and $\{u', v'\} \in A^{(2)}$ such that $\{a_i, b_i\} \to_n \{r_i, s_i\} \Rightarrow_{N, 2^N} \{u', v'\}$ for all $i = 1, \ldots, N$, and $\{u, v\} \Rightarrow_{N, 2^N} \{u', v'\}$. In particular, $\{a_i, b_i\} \Rightarrow_{n+N, 2^N} \{u', v'\}$ for all $i$.*

*Proof.* For each $i = 1, \ldots, N$ choose a sequence $S_i = (u_0^i, u_1^i, \ldots, u_{\lambda(i)}^i)$ from $u$ to $v$ so that $\{a_i, b_i\} \to_n \{u_j^i, u_{j+1}^i\}$ for all $i = 1, \ldots, N$ and all $j < \lambda(i)$. Apply the previous lemma to the sequences $S_1, \ldots, S_N$. □

For integers $n, k \geq 0$ let $C(n, k)$ denote the binomial coefficient. For fixed $m \geq 2$ define $L = L(m)$, $M = M(m)$, $N = N(m)$, $d = d(m)$ and $\ell = \ell(m)$ by

$$\begin{aligned}
L &= C(m+1, 2), \\
M &= C(2m, m) - 1, \\
N &= 2C(M+1, 2), \\
d &= N + 3, \\
\ell &= dM + 2.
\end{aligned}$$

Given a language, let $\Phi_m$ be a formal infinitary sentence in that language (first-order if the language is finite) which asserts

$$\exists x_0 \cdots x_m y z \left[ y \neq z \ \& \ \bigwedge_{0 \leq i < j \leq m} \{x_i, x_j\} \Rightarrow_{2^{6m}, 2^L} \{y, z\} \right].$$

**Lemma 3.4.** *Suppose $\mathcal{V}$ is a class of algebras satisfying the hypotheses of Lemma 3.1. Then the following is true: if $\mathbf{A} \in \mathcal{V}$ and $2 \leq m < \omega$, then either*

1. $\mathbf{A} \models \Phi_m$, *or*
2. *for all $\{a, b\}, \{c, d\} \in A^{(2)}$, $\mathrm{Cg}(a, b) \cap \mathrm{Cg}(c, d) \neq 0_A$ if and only if there exists $\{u, v\} \in A^{(2)}$ such that $\{a, b\} \Rightarrow_{\ell, 4} \{u, v\}$ and $\{c, d\} \Rightarrow_{\ell, 4} \{u, v\}$, where $\ell = \ell(m)$ is as defined above.*

*Proof.* Assume $\mathbf{A} \in \mathcal{V}$ and $\mathbf{A} \not\models \Phi_m$; we shall prove the equivalence in item 2. ($\Leftarrow$) is clear. To prove the other direction, let $a, b, c, d \in A$ with $\mathrm{Cg}(a,b) \cap \mathrm{Cg}(c,d) \neq 0_A$. In fact, we will prove the existence of $\{r,s\}, \{r',s'\}, \{u,v\} \in A^{(2)}$ such that

$$\{a,b\} \to_{dM} \{r,s\} \Rightarrow_{2,4} \{u,v\},$$
$$\{c,d\} \to_{dM} \{r',s'\} \Rightarrow_{2,4} \{u,v\}.$$

To this end, if $n > 0$ and $\{x,y\}, \{z,w\} \in A^{(2)}$ let us say that $\{x,y\}$ and $\{z,w\}$ are $(n,2,4)$-*bounded* (in $\mathbf{A}$) if there exist $\{r,s\}, \{r',s'\}, \{u,v\} \in A^{(2)}$ such that

$$\{x,y\} \to_n \{r,s\} \Rightarrow_{2,4} \{u,v\},$$
$$\{z,w\} \to_n \{r',s'\} \Rightarrow_{2,4} \{u,v\}.$$

Suppose that $\{a,b\}, \{c,d\}$ are *not* $(dM,2,4)$-bounded. By Corollary 3.3, $\{a,b\}, \{c,d\}$ are $(n,2,4)$-bounded for some $n > 0$. Choose $n$ so that $\{a,b\}, \{c,d\}$ are $(n,2,4)$-bounded but are not $(k,2,4)$-bounded for any $k < n$. Thus $n > dM$. Also choose $\{r,s\}, \{r',s'\}, \{u,v\} \in A^{(2)}$ witnessing $(n,2,4)$-boundedness. The remainder of the argument is almost identical to Baker's argument in [1, §8 Case 2]; the reader is urged to read this first.

Let $t = n - dM$. As $\{a,b\} \to_n \{r,s\}$ and $\{c,d\} \to_n \{r',s'\}$ we can choose $\{a_0, b_0\}, \{a'_0, b'_0\} \in A^{(2)}$ such that

$$\{a,b\} \to_t \{a_0, b_0\} \to_{dM} \{r,s\},$$
$$\{c,d\} \to_t \{a'_0, b'_0\} \to_{dM} \{r',s'\}.$$

Thus $\{a_0, b_0\}$ and $\{a'_0, b'_0\}$ are $(dM,2,4)$-bounded but are not $(k,2,4)$-bounded for any $k < dM$. Choose $\{a_i, b_i\}, \{a'_i, b'_i\} \in A^{(2)}$ ($1 \leq i \leq M$) such that

$$\{a_0, b_0\} \to_d \{a_1, b_1\} \to_d \cdots \to_d \{a_{M-1}, b_{M-1}\} \to_d \{a_M, b_M\} = \{r,s\}$$

and similarly for the primed pairs. What is important is that

$$\{a_0, b_0\} \to_{dj} \{a_j, b_j\} \to_{d(M-j)} \{r,s\}$$

for $j = 1, \ldots, M$, and similarly for the primed pairs. For each $j > 0$ choose $d(M-j)$-translations $f_j, f'_j$ which witness $\{a_j, b_j\} \to_{d(M-j)} \{r,s\}$ and $\{a'_j, b'_j\} \to_{d(M-j)} \{r',s'\}$ respectively. We can assume that the notation for the pairs has been coordinated so that $f_j(a_j) = r$ and $f_j(b_j) = s$, and similarly for the primed pairs.

Recall that $\{r,s\} \Rightarrow_2 \{u,v\}$. Choose $k < \omega$, elements $u_0, \ldots, u_k \in A$ and 2-translations $g_0, \ldots, g_{k-1}$ such that $u_0 = u$, $u_k = v$, and $\{u_h, u_{h+1}\} = \{g_h(r), g_h(s)\}$ for $h = 0, \ldots, k-1$.

For $0 \leq i < j \leq M$ and $0 \leq h < k$ let $R_{ij}$ be the sequence from $r$ to $s$ obtained by applying $f_j$ to $(a_j, a_i, b_i, b_j)$, and let $S_{ijh}$ be the sequence from $u_h$ to $u_{h+1}$ obtained by applying $g_h$ to $R_{ij}$ or its reverse. Let $S_{ij}$ be the sequence obtained by concatenating the sequences $S_{ij0}, S_{ij1}, \ldots, S_{ijk-1}$. Thus $S_{ij}$ is a sequence from $u$ to $v$ such that for each adjacent pair $\{x,y\}$ with $x \neq y$, one of the following holds:

1. $\{a_i, a_j\} \to_{d(M-j)+2} \{x, y\}$.
2. $\{b_i, b_j\} \to_{d(M-j)+2} \{x, y\}$.
3. $\{a_i, b_i\} \to_{d(M-j)+2} \{x, y\}$.

In a similar fashion, for $0 \leq i < j \leq M$ we can obtain a sequence $S'_{ij}$ from $u$ to $v$ such that for each adjacent pair $\{x', y'\}$ with $x' \neq y'$, one of the following holds:

4. $\{a'_i, a'_j\} \to_{d(M-j)+2} \{x', y'\}$.
5. $\{b'_i, b'_j\} \to_{d(M-j)+2} \{x', y'\}$.
6. $\{a'_i, b'_i\} \to_{d(M-j)+2} \{x', y'\}$.

In all we get $2C(M+1, 2) = N$ sequences from $u$ to $v$. By Lemma 3.2 there exist $\{u_1, v_1\} \in A^{(2)}$ and 'key links' $\{x_{ij}, y_{ij}\}$ for $S_{ij}$ and $\{x'_{ij}, y'_{ij}\}$ for $S'_{ij}$, $0 \leq i < j \leq M$, such that $\{x_{ij}, y_{ij}\} \Rightarrow_{N, 2^N} \{u_1, v_1\}$ and similarly for the primed key links. Let $\{u_1, v_1\}$ and the key links be chosen and fixed.

Suppose for some $S_{ij}$ that the chosen key link $\{x_{ij}, y_{ij}\}$ satisfies case 3 above, i.e., that $\{a_i, b_i\} \to_{d(M-j)+2} \{x_{ij}, y_{ij}\}$. Then

$$\{a_0, b_0\} \to_{di} \{a_i, b_i\} \to_{d(M-j)+2} \{x_{ij}, y_{ij}\} \Rightarrow_N \{u_1, v_1\},$$

and hence $\{a_0, b_0\} \Rightarrow_{dM-1} \{u_1, v_1\}$ as $di + d(M-j) + 2 + N \leq dM - 1$. Likewise, if for some $S'_{i'j'}$ the chosen key link satisfies case 6 above, then $\{a'_0, b'_0\} \Rightarrow_{dM-1} \{u_1, v_1\}$. Thus if cases 3 and 6 both occur, then $\{a_0, b_0\}$ and $\{a'_0, b'_0\}$ would be $(dM-1, 2, 4)$-bounded by Corollary 3.3. As this is not the case, either case 3 or case 6 never occurs. By symmetry we can assume that case 3 never occurs. Thus for all $0 \leq i < j \leq M$, the chosen key link of $S_{ij}$ satisfies case 1 or 2. As $M + 1 = C(2m, m)$, which is at least as large as the Ramsey number for $(m+1)$-element monochromatic subsets of 2-coloured complete graphs, Ramsey's theorem tells us there exists an $(m+1)$-element subset $J \subseteq \{0, 1, \ldots, M\}$ such that either all sequences $S_{ij}$ with $i, j \in J$ and $i < j$ have their key links satisfying case 1, or all have their key links satisfying case 2. Again for concreteness let us assume that all have their key links satisfying case 1. That is, if $i, j \in J$ with $i < j$ then

$$\{a_i, a_j\} \to_{d(M-j)+2} \{x_{ij}, y_{ij}\} \Rightarrow_{N, 2^N} \{u_1, v_1\}$$

and therefore $\{a_i, a_j\} \Rightarrow_{dM} \{u_1, v_1\}$. Then by Corollary 3.3 there exists $\{u', v'\} \in A^{(2)}$ such that $\{a_i, a_j\} \Rightarrow_{dM+L, 2^L} \{u', v'\}$ for all $i, j \in J$ with $i < j$. As in [1, p. 228], the following inequalities

$$M + 1 \leq \frac{12}{11}\left(\frac{2^{2m}}{(\pi m)^{1/2}}\right),$$
$$d \leq (M+1)^2,$$
$$\lceil \log_2 L \rceil \leq m,$$
$$\left(\frac{12}{11}\right)^3 \leq \frac{\pi^{3/2}}{4},$$

imply $\ell\lceil\log_2 L\rceil \leq (M+1)^3 m \leq 2^{6m-2}$ and so certainly $dM + L \leq 2^{6m}$. But this means $\mathbf{A} \models \Phi_m$, contradicting our initial assumption. This proves that $\{a,b\}, \{c,d\}$ are $(dM, 2, 4)$-bounded whenever $\mathrm{Cg}(a,b) \cap \mathrm{Cg}(c,d) \neq 0_A$, proving the lemma. □

**Corollary 3.5.** *Let $\mathcal{V}$ be a class of algebras satisfying the hypotheses of Lemma 3.1. For all $m \geq 2$, if $\mathbf{A} \in \mathcal{V}$ is subdirectly irreducible and $|A| > m$, then $\mathbf{A} \models \Phi_m$.*

*Proof.* Suppose $\mathbf{A} \not\models \Phi_m$. Then by the previous lemma and because $\mathbf{A}$ is subdirectly irreducible, for all $\{a,b\}, \{c,d\} \in A^{(2)}$ there exists $\{u,v\} \in A^{(2)}$ such that $\{a,b\} \Rightarrow_\ell \{u,v\}$ and $\{c,d\} \Rightarrow_\ell \{u,v\}$. It follows by induction on $k$ that if $S \subseteq A^{(2)}$ with $|S| \leq 2^k$ then there exists $\{u,v\} \in A^{(2)}$ such that $\{a,b\} \Rightarrow_{\ell k} \{u,v\}$ for all $\{a,b\} \in S$. Now choose $U \subseteq A$ with $|U| = m+1$ and apply this observation to the set $S$ of all two-element subsets of $U$. Since $|S| = L \leq 2^{\lceil \log_2 L \rceil}$ and $\ell\lceil\log_2 L\rceil \leq 2^{6m-2}$ (see the proof of the previous lemma) we get $\{u,v\} \in A^{(2)}$ such that $\{a,b\} \Rightarrow_{2^{6m-2}} \{u,v\}$ for all distinct $a, b \in U$. By Corollary 3.3, there exists $\{u', v'\} \in A^{(2)}$ such that $\{a,b\} \Rightarrow_{2^{6m}, 2L} \{u', v'\}$ for all distinct $a, b \in U$. This proves $\mathbf{A} \models \Phi_m$ after all. □

## 4. The finite basis theorem

In an arbitrary algebra $\mathbf{A}$ we let $M(x, y, z, w)$ denote the 4-ary relation on $A$ defined by $\mathrm{Cg}(x, y) \cap \mathrm{Cg}(z, w) \neq 0_A$. $\mathbf{A}$ is *finitely subdirectly irreducible* if it satisfies $\forall xyzw[(x \neq y \ \& \ z \neq w) \to M(x, y, z, w)]$. If $\mathcal{K}$ is a class of algebras, then $\mathbb{S}_{\mathbb{I}}(\mathcal{K})$ denotes the class of all subdirectly irreducible members of $\mathcal{K}$ while $\mathbb{F}_{\mathbb{S}\mathbb{I}}(\mathcal{K})$ denotes the class of all finitely subdirectly irreducible members of $\mathcal{K}$.

The following lemma is due to B. Jónsson (see [7, Theorem 1] or [8, Lemma 7.2]).

**Lemma 4.1.** *Suppose $\mathcal{V}$ is a variety of algebras and $\mathcal{H}$ is a strictly elementary class containing $\mathcal{V}$. If there exists an elementary class $\mathcal{K}$ such that $\mathbb{S}_{\mathbb{I}}(\mathcal{H}) \subseteq \mathcal{K}$ and $\mathcal{V} \cap \mathcal{K}$ is strictly elementary, then $\mathcal{V}$ is finitely based.*

The next lemma is presumably folklore; we show how to deduce it from Lemma 4.1.

**Lemma 4.2.** *Suppose $\mathcal{V}$ is a variety in a finite language and suppose there exists $m < \omega$ such every $\mathbf{A} \in \mathbb{S}_{\mathbb{I}}(\mathcal{V})$ has cardinality less than $m$. If there exists a strictly elementary class $\mathcal{H}_0$ such that (i) $\mathcal{V} \subseteq \mathcal{H}_0$ and (ii) the relation $M(x, y, z, w)$ is definable throughout $\mathcal{H}_0$ by a first-order formula, then $\mathcal{V}$ is finitely based.*

*Proof.* First observe that $\mathbb{F}_{\mathbb{S}\mathbb{I}}(\mathcal{V}) = \mathbb{S}_{\mathbb{I}}(\mathcal{V})$. Indeed, if $\mathbf{A} \in \mathbb{F}_{\mathbb{S}\mathbb{I}}(\mathcal{V})$ is finite then $\mathbf{A} \in \mathbb{S}_{\mathbb{I}}(\mathcal{V})$ automatically. Assume that there exists an infinite $\mathbf{A} \in \mathbb{F}_{\mathbb{S}\mathbb{I}}(\mathcal{V})$; choose $S \subseteq A$ such that $|S| = m + 1$ where $m$ is the maximum cardinality of the members of $\mathbb{S}_{\mathbb{I}}(\mathcal{V})$. Because $\mathbf{A}$ is finitely subdirectly irreducible, there exist distinct $c, d \in A$ such that $(c, d) \in \mathrm{Cg}(x, y)$ for all distinct $x, y \in S$. Thus if $\theta$ is a congruence of $\mathbf{A}$ maximal with respect to omitting $(c, d)$, then $\mathbf{A}/\theta$ is a member of $\mathbb{S}_{\mathbb{I}}(\mathcal{V})$ having cardinality greater than $m$, a contradiction.

Now let $\mathcal{K} = \mathbb{F}_{\mathbb{SI}}(\mathcal{V})$ and $\mathcal{H} = (\mathcal{H}_0 \setminus \mathbb{F}_{\mathbb{SI}}(\mathcal{H}_0)) \cup \mathcal{K}$. $\mathbb{F}_{\mathbb{SI}}(\mathcal{H}_0)$ and $\mathcal{K}$ are strictly elementary because of the definability of $M(x,y,z,w)$ in $\mathcal{H}_0$ and the bounded size of members of $\mathcal{K}$ respectively; therefore $\mathcal{H}$ is strictly elementary. It follows from Lemma 4.1 that $\mathcal{V}$ is finitely based. □

**Theorem 4.3.** *Suppose $\mathcal{V}$ is a congruence meet-semidistributive variety in a finite language. If there exists $m < \omega$ such that every $\mathbf{A} \in \mathbb{S}_{\mathbb{I}}(\mathcal{V})$ has cardinality less than $m$, then $\mathcal{V}$ is finitely based.*

*Proof.* Let $m < \omega$ be larger than $|A|$ for every $\mathbf{A} \in \mathbb{S}_{\mathbb{I}}(\mathcal{V})$. By Theorem 2.1 we can choose a finite family $\{\langle s_p, t_p \rangle : p \in T\}$ of pairs of ternary terms witnessing Theorem 2.1(6). Let $\mathsf{s}_p, \mathsf{t}_p$ ($p \in T$) be new 3-ary operations symbols, let $\mathcal{L}'$ be the extension of the language of $\mathcal{V}$ to these new symbols, and let $\mathcal{V}'$ be the obvious expansion of $\mathcal{V}$ to $\mathcal{L}'$ defined by $\mathsf{s}_p(x,y,z) \approx s_p(x,y,z)$ and $\mathsf{t}_p(x,y,z) \approx t_p(x,y,z)$ ($p \in T$). Thus $\mathcal{V}'$ is a variety in a finite language and is term-equivalent to $\mathcal{V}$. We shall work with $\mathcal{V}'$ in place of $\mathcal{V}$. Let $\mathcal{V}^*$ be the (strictly elementary) class in the language $\mathcal{L}'$ defined by the conditions in Theorem 2.1(6) with respect to the family $\{\langle \mathsf{s}_p, \mathsf{t}_p \rangle : p \in T\}$; thus $\mathcal{V}' \subseteq \mathcal{V}^*$ and $\mathcal{V}^*$ satisfies the hypotheses of Lemma 3.1.

Recall the definitions of $\ell$ and $\Phi_m$ preceding Theorem 3.4. Because $\mathcal{L}'$ is finite, $\Phi_m$ is first-order; likewise, there is a first-order formula $\mu(x,y,z,w)$ which asserts in any $\mathbf{A} \in \mathcal{V}^*$ that $x \neq y$ and $z \neq w$ and there exists $u \neq v$ such that $\{x,y\} \Rightarrow_{\ell,4} \{u,v\}$ and $\{z,w\} \Rightarrow_{\ell,4} \{u,v\}$. Any model of $\Phi_m$ has a subdirectly irreducible homomorphic image of cardinality greater than $m$; hence $\mathcal{V}' \models \neg\Phi_m$. Let $\mathcal{H}_0 = \{\mathbf{A} \in \mathcal{V}^* : \mathbf{A} \models \neg\Phi_m\}$. $\mathcal{H}_0$ is strictly elementary and contains $\mathcal{V}'$. By Lemma 3.4, the relation $M(x,y,z,w)$ is defined in $\mathcal{H}_0$ by the formula $\mu$. Thus $\mathcal{V}'$ (and therefore also $\mathcal{V}$) is finitely based by Lemma 4.2. □

## 5. The finite spectrum of $\mathbb{S}_{\mathbb{I}}(\mathcal{V})$

In this section we show that if $\mathcal{V}$ is a locally finite congruence meet-semidistributive variety in a finite language $\mathcal{L}$, then the gaps in the set $\{|A| : \mathbf{A} \in \mathbb{S}_{\mathbb{I}}(\mathcal{V}), \mathbf{A} \text{ finite}\}$ are controlled by the free spectrum of $\mathcal{V}$ and the complexity of $\mathcal{L}$. As a consequence, we obtain the algorithm promised in Theorem B of the Abstract.

If $\mathcal{V}$ is such a variety then define $\mathcal{L}'$ and $\mathcal{V}'$ as in the proof of Theorem 4.3. Let $M'$ be the maximum of the arities of the symbols in $\mathcal{L}'$. Let $\mathbf{A}$ be a member of $\mathcal{V}$ and let $\mathbf{A}'$ be its unique expansion in $\mathcal{V}'$. If $\{a,b\} \to_k \{c,d\}$ in $\mathbf{A}'$, then this is witnessed by a $k$-translation of $\mathbf{A}'$; this translation can be defined from a term of $\mathbf{A}'$ (or of $\mathbf{A}$) using at most $k(M'-1)$ parameters. Similarly, if $\{a,b\} \Rightarrow_{k,n} \{c,d\}$ in $\mathbf{A}'$ then this can be witnessed by at most $kn(M'-1)$ parameters.

Now suppose that $\mathbf{A} \in \mathcal{V}$ is subdirectly irreducible and $|A| > m \geq 2$. Then $\mathbf{A}'$ is also subdirectly irreducible and $|A'| = |A| > m$. By Corollary 3.5, $\mathbf{A}' \models \Phi_m$. Thus if $L = C(m+1, 2)$ then there exist $a_0, \ldots, a_m, u, v \in A$ with $u \neq v$ such that $\{a_i, a_j\} \Rightarrow_{2^{6m}, 2^L} \{u, v\}$ for all $0 \leq i < j \leq m$. It follows that there is a subalgebra $\mathbf{B}'$ of $\mathbf{A}'$ generated by $a_0, \ldots, a_m$ and at most $L \cdot 2^{6m+L}(M'-1)$ other parameters such that $\mathbf{B}' \models \Phi_m$. $\mathbf{B}'$ has a homomorphic image which is subdirectly irreducible and has cardinality greater than $m$; the reduct of $\mathbf{B}'$ to the language of $\mathcal{V}$ also has this property. Since $L \cdot 2^{6m+L}(M'-1) + m + 1 \leq 2^{m^2+7m}M'$ this proves:

**Theorem 5.1.** *Let $\mathcal{V}$ be a locally finite congruence meet-semidistributive variety in a finite language. Let $M$ be the maximum of the arities of the fundamental operations of $\mathcal{V}$ and let $M' = \max(M, 3)$. If $2 \leq m < \omega$ and there exists a member of $\mathbb{S}_\mathbb{I}(\mathcal{V})$ of cardinality greater than $m$, then there exists a finite member $\mathbf{B} \in \mathbb{S}_\mathbb{I}(\mathcal{V})$ such that $m < |B| \leq f_\mathcal{V}(2^{p(m)}M')$, where $p(m) = m^2 + 7m$ and $f_\mathcal{V}(n)$ denotes the cardinality of the $\mathcal{V}$-free algebra on $n$ generators.*

By a *concrete language* we mean the graph $L$ of a function from a finite initial segment of $\omega$ into $\omega$, where each $(i, k) \in L$ is construed as an operation symbol of arity $k$. A *concrete algebra* is an algebra whose language is concrete and whose universe is a finite initial segment of $\omega$. The set $\mathbb{C}\text{ONC}$ of all concrete algebras may be effectively encoded as a recursive subset of $\omega$. Thus we may speak of subsets of $\mathbb{C}\text{ONC}$ as being recursive or r.e.

A variety $\mathcal{V}$ is *residually less than $m$* if $|B| < m$ for all $\mathbf{B} \in \mathbb{S}_\mathbb{I}(\mathcal{V})$. Define

$$\begin{aligned}
\mathcal{S} &= \{\mathbf{A} \in \mathbb{C}\text{ONC} : \mathbf{HSP}(\mathbf{A}) \text{ is congruence meet-semidistributive}\} \\
\mathcal{S}_m &= \{\mathbf{A} \in \mathcal{S} : \mathbf{HSP}(\mathbf{A}) \text{ is residually less than } m\}, \quad (1 < m < \omega) \\
\mathcal{S}_\omega &= \bigcup_{m < \omega} \mathcal{S}_m.
\end{aligned}$$

**Corollary 5.2.** *The set $\{(m, \mathbf{A}) : 1 < m < \omega \text{ and } \mathbf{A} \in \mathcal{S}_m\}$ is recursive. $\mathcal{S}_\omega$ is r.e. Hence among finite algebras $\mathbf{A}$ in a finite language which generate congruence meet-semidistributive varieties, the property "$\mathbf{HSP}(\mathbf{A})$ has a finite residual bound" is semi-decidable.*

*Proof.* $\mathcal{S}$ is recursive by Theorem 2.1($1 \Leftrightarrow 3$). Given $m < \omega$ and $\mathbf{A} \in \mathcal{S}$, put $k = |A|$, let $M'$ and $p(m)$ be as in the previous theorem, let $L$ be the language of $\mathbf{A}$, and define $T = k^{k^{(2^{p(m)}M')}}$. Thus $f_\mathcal{V}(2^{p(m)}M') \leq T$, so if there exists $\mathbf{B} \in \mathbb{S}_\mathbb{I}(\mathbf{HSP}(\mathbf{A}))$ with $|B| \geq m$ then such $\mathbf{B}$ exists with $|B| \leq T$. This can be effectively tested by enumerating the finitely many concrete algebras $\mathbf{B}_1, \ldots, \mathbf{B}_N$ of type $L$ and such that $m \leq |B_i| \leq T$, and then individually testing each $\mathbf{B}_i$ as to whether it (i) is subdirectly irreducible and (ii) is in $\mathbf{HSP}(\mathbf{A})$. Item (ii) can be determined effectively by an old argument of J. Kalicki [9]. This proves the first claim. The second claim is a consequence of the first. □

## 6. Conclusion

Looking forward, it is natural to ask for a finite basis theorem which incorporates both McKenzie's theorem for congruence modular varieties and our Theorem 4.3 for congruence meet-semidistributive varieties. A possible setting for such a generalization is provided by tame congruence theory. Our theorem handles the (finitely generated) varieties which omit types **1** and **2**, while McKenzie's theorem, together with [5, Theorem 10.4], handles the varieties which omit types **1** and **5**.

**Problem 2.** Is Problem 1 true for finitely generated varieties which omit type **1**?

In an algebra let $C(x, y, z, w)$ denote the 4-ary relation "$[\mathrm{Cg}(x,y), \mathrm{Cg}(z,w)] \ne 0$" where $[\_, \_]$ denotes the TC commutator. $C(x, y, z, w)$ is identical to $M(x, y, z, w)$ in any congruence meet-semidistributive variety, but this is no longer true in varieties that admit abelian phenomena. A key step in McKenzie's proof of his finite basis theorem was the demonstration that $C(x, y, z, w)$ is definable in any congruence modular variety in a finite language which has a finite residual bound.

**Problem 3.** Suppose $\mathcal{V}$ is a variety in a finite language which omits type **1** and has a finite residual bound. Is the relation $C(x, y, z, w)$ definable in $\mathcal{V}$ by a first-order formula?

As soon as Baker revealed his finite basis theorem, general algebraists sought and found proofs which were simpler in that they avoided the Ramsey argument. The recent proof of Baker and Ju Wang [2] is perhaps the nicest example.

**Problem 4.** Find an elementary proof of Theorem 4.3 which does not require a combinatorial analysis of congruence generation.

DEPARTMENT OF PURE MATHEMATICS, UNIVERSITY OF WATERLOO, WATERLOO, ONTARIO N2L 3G1, CANADA
*E-mail address*: `rdwillar@gillian.math.uwaterloo.ca`